Andrew Aberdein

## CLASSICAL RECAPTURE

The recapture relationship is an important element to any understanding of the connexion between different systems of logic. Loosely speaking, one system of logic recaptures another if it is possible to specify a subsystem of the former system which exhibits the same patterns of inference as the latter system.[1] In particular if a relationship of this kind can be shown to exist between a non-classical logic and classical logic (henceforth **K**), the non-classical system is said to exhibit classical recapture. This has been invoked by several proponents of non-classical logics to argue that their system retains **K** as a limit case, and is therefore a methodologically progressive successor to **K**. In this paper I shall advance and defend a new and more precise account of recapture and the character of its reception by the proponents of the recapturing system. I shall then indicate some of the applications of classical recapture which this account makes possible.

Logics can be presented in many different ways—natural deduction presentations, sequent calculi, various axiom systems, and so forth—but three basic types of presentation may be distinguished: logistic systems, which codify logical truths; consequence systems, which codify valid arguments; and deductive systems, which codify proofs.[2] My concern is with substantive divergence amongst logical systems intended for the formalization of rational argumentation. Although logistic systems may be adequate for some purposes, such as codifying the truths of arithmetic, they are too coarse-grained to capture all the differences with which I am concerned.[3] Conversely, deductive systems offer too fine-grained a classification: differences which occur only at this level are outside the scope of my inquiry. Therefore my attention may be safely restricted to consequence systems.

---

[1] The earliest usage I can find of the word 'recapture' to describe a relationship of this kind is Priest 1987 p146, although such relationships have been discussed in other terms for much longer. Sometimes this has been in a weaker sense, as the reproduction of the theorems of the prior system, or in a stronger sense, as the reproduction of the proofs of that system.

[2] Corcoran 1969 pp154ff. Equivalently, the gross and delicate proof theory of Tennant 1996 pp351f consist in the analysis of consequence and deductive systems respectively.

[3] For example, **K** has the same theorems as the relevant system **R**.

Yet there is more to the formalization of argumentation than the presentation of a formal system. We must also be concerned with the parsing theory, by which translation to and from the formal system is effected, and various background assumptions.[4] Taken together with the formal system, these factors constitute a logical theory, by which the system may be promoted. Sequences of logical theories may be considered as logical research programmes, characterized by the retention of an irrevisable hard core.[5] Programmes have heuristic as well as theoretical content: methods for constructing more successful theories while protecting the hard core from pressure for revision.

In Aberdein 1998 I introduced and defended an account of the equivalence of consequence systems. This utilised a schematized representation of such systems, $L_i$, as couples, $<W_i, V_i>$, where $W_i$ is the class of well-formed formulæ of the language underpinning logic $L_i$ and $V_i$ is the class of valid inferences of $L_i$ (a subclass of the class of sequents defined on $W_i$). Equivalence consists in a pair of surjections between the classes of wffs of the systems which preserve the partitions of the classes of inferences into valid and invalid subclasses. I proceeded to define a means of contracting a formal system:

DEFINITION 1: $L_1$ is a *proper fragment* of $L_2$ iff $L_1$ and $L_2$ are inequivalent, $W_1$ is defined on a proper subset of the class of constants of $L_2$ and $V_1$ contains precisely those elements of $V_2$ which contain only elements of $W_1$.

Hence, fragmentation is the inverse of conservative extension. However, fragments are not the only sort of contractions that may be defined upon formal systems; the definition may be generalized as follows:

DEFINITION 2: $L_1$ is a *proper subsystem* of $L_2$ iff $L_1$ and $L_2$ are inequivalent, $W_1$ is a proper subset of $W_2$ and $V_1$ contains precisely those elements of $V_2$ which contain only elements of $W_1$.

The metaphors of strength, size and inclusion which so often illustrate the mereology of logic suffer from an ambiguity: there is a tension between a deductive characterization, a measure of how much may be deduced from how little, and an expressive characterization, a measure of the subtlety of the distinctions which can be preserved.[6] An increase in one may represent a decrease in the other. Hence, 'subsystem of $L$' has sometimes been used to designate a system axiomatized by a subset of the axioms of $L$, or with a deducibility relation which is a subrelation of that of $L$. The definition of subsystem adopted above reverses this usage, making explicit the generalization of the definition of fragment, but rendering these 'subsystems' supersystems, the inverse of subsystems. In short, fragments are exclusively generated by reducing the set of constants upon which the class of wffs is based, but subsystems may also be generated by reducing the class of wffs in some other way. For example, **K** is a subsystem of intuitionistic logic (henceforth **J**): for **K** may be thought of as the system resulting from **J** by the exclusion of all undecidable formulæ, as may be achieved by adding the law of excluded middle to the axioms of **J**, or the rule of double-negation elimination to the definition of its deducibility relation.

---

[4] See Thagard 1982 and Resnik 1985 for contrasting accounts of this material.

[5] I discuss at length the application to logic of these Lakatosian ideas in Aberdein 1999.

[6] *Cf*. the distinction between expressive power and deductive power drawn by Rautenberg (1987 p.*xvi*), discussed in Beziau 1997 pp5f.

This apparatus provides the means for a formal account of recapture.

DEFINITION 3: $L_1$ *recaptures* $L_2$ iff there is a proper subsystem of $L_1$, $L_1$\*, which is defined in terms of a constraint on $W_1$ finitely expressible in $L_1$, and which is equivalent to $L_2$. If $L_2$ is **K**, then $L_1$ is a *classical recapture logic*.

Which is to say that if one system recaptures another we may express within it some finite constraint by which a subsystem equivalent to the recaptured system may be generated. For example, we can see that **J** is a classical recapture logic, with the constraint of decidability. The relevant system **R** and quantum logic also recapture **K**, with constraints of negation consistency and primality, and compatibility, respectively. Indeed, many non-classical logics are classical recapture logics: exactly which will turn on which constraints are deemed expressible. It has even been suggested that the recapture of **K** is a necessary criterion of logicality, in which case all logics would be classical recapture logics.[7]

Some non-classical logicians embrace classical recapture; others attempt to reject it; while others see recapture results as motivating the reduction of the recapturing system to a conservative extension. Thus, before recapture can contribute to the kinematics of logic, we must distinguish amongst the variety of responses that advocates of a logical system may make to the status of their system as a (classical) recapture logic. I shall order these responses by analogy with a spectrum of political attitudes: radical left, centre left, centre right and reactionary right. This is a formal not a sociological analogy: I do not intend to imply that views on logic may be correlated to political allegiance (*pace* some sociologists of scientific knowledge). The most extreme of these attitudes is the radical left: formal repudiation of recapture status. Individuals of this tendency deny that their system recaptures the prior system, claiming that no suitable recapture constraint is expressible in the new system. If classical recapture were a criterion of logicality, then a radical left response could only be embraced by quitting the discipline of logic. Yet such a criterion must be open to doubt, since some familiar programmes include proponents from the radical left. For example, Nuel Belnap and Michael Dunn's argument that relevant logic does not recapture **K** places their relevantist in this camp.[8] The subordination of logic to mathematics by some proponents of **J** may also be understood as preventing classical recapture.

---

[7] 'Perhaps … any genuine 'logical system' should contain classical logic as a special case' van Benthem 1994 p135.

The less radical centre left acknowledge the formal satisfaction of recapture, but deny its significance. Proponents of this stance argue that the formal equivalence between a subsystem of their system and another system is irrelevant, since the other system cannot be understood as formalizing anything intelligible in terms of their theory. Hence some advocates of **J** regard the double-negation translation of **K** into their system as no more than a curiosity, since they reject the cogency of classical concepts.[9] Whereas the radical left presume a logical incompatibility between the recapture result and indispensible formal components of the research programme, the centre left claim an heuristic incompatibility with indispensible non-formal components of the research programme. To defend a position on the centre left one must demonstrate that conceding more than a technical significance to recapture will induce an intolerable tension between successful problem-solving within the programme and the retention of its key non-formal components, such as the central aspects of its parsing theory.

On the centre right recapture is embraced as evidence of the status of the new system as a methodologically progressive successor. The meaning invariance of all key terms is welcomed in this context, and recapture is understood as establishing the old system as a limit case of its successor. The centre right hold with Einstein that '[t]here could be no fairer destiny for any … theory than that it should point the way to a more comprehensive theory in which it lives on, as a limiting case'.[10] Most non-classical logics have been defended in these terms by at least some of their advocates: for example, Hilary Putnam's quondam advocacy of quantum logic was of this character, as is Graham Priest's support for paraconsistent logic.[11]

Least radical of all are the reactionary right, who argue that the subsystem of the new system equivalent to the old system is actually a proper fragment of the new system, that is that the new system should be understood as extending the old system. Hence the *status quo* is maintained: the old system is still generally sound, but can be extended to cover special cases. Many ostensibly non-classical programmes have at some stage been promoted as conservative extensions of **K**: for example, Maria Louisa Dalla Chiara's modal quantum logic **B$^O$** or Robert Meyer's classical relevant system **R$^¬$**.[12] Modal logic may be understood as having successfully completed a move from the centre right to the reactionary right: although it is now understood as extending **K**, its early protagonists conceived it as a prospective successor system.[13]

---

[8] Anderson, Belnap & Dunn 1992 §80.4.5 p505. In this case the situation is complicated by their claim not to embrace the radical left stance themselves; rather they attribute it to a position they wish to criticize.

[9] '[I]ntuitionists … deny that the [classical] use [of the logical constants] is coherent at all' Dummett 1973b p398. But see Dummett 1973a p238 for a more conciliatory intuitionist response to recapture.

[10] A. Einstein 1916 *Relativity: The special and general theory* quoted in Popper 1963 p32.

[11] Putnam 1968 p184; Priest 1987 pp146ff.

Different logical research programmes encompass different political complexions: some are clearly associated with one stance, whether for technical or historical reasons, in others there is dispute as to which approach is appropriate. Two further points may serve to reinforce the political analogy: programmes appear to drift to the right as they grow older, and there is a strong community of interest between the two ends of the spectrum. The reactionary agrees with the left-wingers that the constants of the new system have different meanings from those of the old. The difference is that the left wing think that the new meanings must replace the old, whereas reactionaries believe that they can be assimilated into an augmented system through employment alongside the old meanings. The greater the difference between the new and the old constants, the more difficult it is to maintain a centrist position.

The full range of options may be seen more clearly as a flow chart:

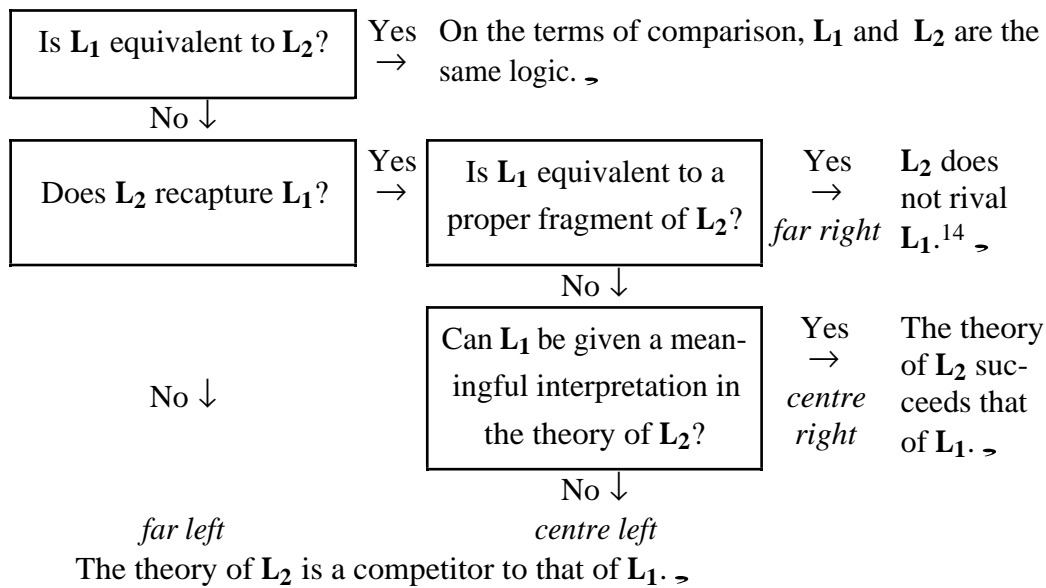

This chart has been devised to display the consequences of a change of theory in which a specific formal system ($L_2$) replaces another ($L_1$). However, it should be stressed that, in the practical development of logical research programmes, a dialectic exists between the choice of formal system and the attitude taken to the recapture of the prior system. Hence, providing that enough of the formal system remains within the revisable part of a logical research programme, there are always two alternatives: embrace the consequences of the formal system, or change the system to resist them.

With this picture in place, I can begin to outline some of the uses to which it may be put. In the first place, we now have the resources to draw some fundamental distinctions between different sorts of theory change. An important feature of the flow chart is that its first three questions can be answered purely by comparison of the formal systems $L_1$ and $L_2$, but the fourth question, 'Can $L_1$ be given a meaningful interpretation in the theory of $L_2$?', requires an appeal to the theories by which the systems are advanced, and perhaps the research programme behind that. Hence, while certain outcomes are necessitated by formal features, other outcomes are underdetermined by such data alone. Solely on formal data we can observe that rivalry must occur unless one system conservatively extends the other, and that competition must occur unless one system recaptures the other. However, broader consideration is required if more than these weak sufficiency conditions for the supplement/rival and successor/competitor distinctions are sought. Indeed, logical theories can be rivals even when the embedded systems are related by conservative extension, or even equivalence: for example, $R^\neg$ conservatively extends $K$, but its promotion would presume a radically non-classical parsing theory, and many systems of logic have more than one alternative semantics.[15] Yet, presuming that the remainder of the theory changes no more than necessary, a clear taxonomy of the consequences of different species of logical revision may be seen to emerge.

Secondly, recapture may be assimilated to accounts of theory change from the philosophy of science. For example, Arthur Fine has defended an account of semantic stability across theory change which requires a shared characterization of the contested term and specifiable constraints within the new theory in which the relationships of the old theory still obtain.[16] Thus, when augmented with a suitably general schema for characterizing the logical constants,[17] centre right recapture provides a defence of the meaning invariance of logical constants across changes of logical theory. This provides a key weapon against claims either that logical reform is merely the vapid rearrangement of notation or that it can only be achieved by an incommensurable change of subject.

---

[15] Susan Haack acknowledges the possibility of such a limitation to her attempt to define rivalry on purely syntactic grounds, although her choice of examples downplays its likelihood (Haack 1974 p6).

[16] Fine 1967 pp237f.

[17] Such as that advanced in Dosen 1987 pp366ff.

Finally, we now have the means to get to grips with the methodological claims made for recapture as a hallmark of progress. The methodological implications of the recapture result are twofold: if the centre right attitude is defensible, there will be good grounds for arguing that the new system is a progressive successor to the old; but progress may come by other means. If the old system is sufficiently wrong, recapture may just be a means of continuing in an entrenched error; the only progressive attitude would be on the centre left. Conversely, by the reactionary response, recapture may be a means for the old system to overtake the new, and thereby prove the more progressive. So recapture can be a sign of progress, but it cannot be an infallible touchstone.